%bijstern.tex: a Plain TeX file by Shalosh B. Ekhad and Doron Zeilberger
%A Bijective Proof of Richard Stanley's Observation that the sum of the cubes of the n-th row of`
%Stern's  Diatomic array equals 3*7^(n-1)

%begin macros

\baselineskip=14pt
\parskip=10pt

\magnification=\magstephalf
\def\A{{\cal A}}
\def\B{{\cal B}}
\def\C{{\cal C}}
\def\S{{\cal S}}

\def\1{{\overline{1}}}
\def\2{{\overline{2}}}
\parindent=0pt
\overfullrule=0in

\def\frac#1#2{{#1 \over #2}}
%\headline={\rm  \ifodd\pageno  \RightHead  \else  \LeftHead  \fi}
%\def\RightHead{\centerline{
%Title
%}}
%\def\LeftHead{ \centerline{Doron Zeilberger}}
%end macros
\centerline
{\bf
A Bijective Proof of Richard Stanley's Observation that the sum of the cubes 
 }
\centerline
{\bf
of the n-th row of Stern's Diatomic Array equals 3 times 7 to the power n-1
 }

\bigskip
\centerline
{\it Shalosh B. EKHAD and Doron ZEILBERGER}

\bigskip

{\bf Abstract}: In a delightful article, Richard Stanley derived, algebraically, the  ``{\it surprisingly simple}''
formula, $3 \cdot 7^{n-1}$, for the sum of the cubes of the $n$-th row of {\it Stern's diatomic array}.
In this  note, we find an {\it elegant} bijective proof of this surprising fact, that
{\bf explains} it and gives {\bf insight}. The {\bf novelty} is that
this gorgeous bijection was {\bf discovered} by a computer (SBE), with {\it minimal} guidance by a human (DZ).
This {\bf debunks} the {\it conventional wisdom}, held by some {\it human supremacists}, that computers
can only compute, but they can't give insight.

In [S], the following double sequence $a(n,k)$ is defined
$$
\sum_{k=0}^{2(2^n-1)} \, a(n,k) x^k \, = \, \prod_{i=0}^{n-1} \left ( \, 1+x^{2^i} +x^{2 \cdot 2^i} \, \right )  \quad.
$$

For any vector in $a \in \{0,1,2\}^n$, define the {\bf valuation}, $v(a)$, by
$$
v(a)= \sum_{i=1}^{n} a_i \, 2^{i-1} \quad,
$$

then, obviusly, $a(n,k)$ counts the number of vectors in   $\{0,1,2\}^n$ with valuation $k$, and
hence $a(n,k)^r$ counts the number of $r \times n$ matrices with entries in $\{0,1,2\}$ such that the
valuation of {\it each} of its $r$ rows equals $k$. Hence
$$
A_r(n):=\sum_{k=0}^{2(2^n-1)} a(n,k)^r \quad,
$$
counts the {\bf set}, let's call it $\A_r(n)$, of all $r \times n$ matrices with entries in $\{0,1,2\}$  such
that all rows have identical valuations.

Stanley ([S], p.4)  proved (among many other results), the

{\bf Surprisingly Simple Fact}: $A_3(n)\,=\, 3 \cdot 7^{n-1}$ \quad .

This {\bf cries} for a {\bf bijective proof}! The right side obviously counts the set of
words in the alphabet $\{1,2,3,4,5,6,7\}$ of length $n$, whose {\bf last} letter is in $\{1,2,3\}$.
Let's call this set $\S(n)$.

{\bf Theorem:} There exists an elegant bijection, implemented in procedure {\tt StanBij(M)} in the Maple package

{\tt http://www.math.rutgers.edu/\~{}zeilberg/tokhniot/BijectionStern.txt} \quad ,

between the set $\A_3(n)$ and $\S(n)$, whose inverse is given by {\tt InvStanBij(w)}, also implemented there.

The actual definition of the bijection can be easily read directly from the Maple source code, let's just describe how we constructed it.
The same methodology should work to provide bijective proofs to any of the recurrences listed in [S], but the details
will not be as elegant.

The set $\A_3(3)$ has $3 \cdot 7^2=147$ members, but the set obtained by deleting the third column happens to only
have $63$ elements. In fact this is true if you take $A_3(n)$ and remove the last $n-2$ columns, you always get
the same set. Let's call this set of $3 \times 2$ matrices, $\B_3$.

On the other hand the set of $3 \times 1$ matrices obtained by only retaining the first column of the members
of $\A_3(n)$ (for $n\geq 2$) has exactly $9$ elements, let's call it $\C_3$.
$$
\C_3 \, =\,
\left \{ \left ( \matrix{ 0 \cr 0 \cr 0} \right )  \, , \,
 \left ( \matrix{ 0 \cr 0 \cr 2} \right )  \, , \,
 \left ( \matrix{ 0 \cr 2 \cr 0} \right ) \, , \,
 \left ( \matrix{ 0 \cr 2 \cr 2} \right ) \, , \,
 \left ( \matrix{ 1 \cr 1 \cr 1} \right )\, , \,
\left ( \matrix{ 2 \cr 0 \cr 0} \right )\, , \,
 \left ( \matrix{ 2 \cr 0 \cr 2} \right )\, , \,
 \left ( \matrix{ 2 \cr 2 \cr 0} \right )\, , \,
 \left ( \matrix{ 2 \cr 2 \cr 2} \right ) \right \} \quad .
$$

Given a matrix $M$ in $\A_r(n)$, we define $v(M)$ to be the vector of length $r$ consisting of the
valuations of its rows. If the smallest entry is larger than zero, we subtract it from all the entries,
getting a {\bf reduced valuation vector}, so one of its entries must be $0$.

It so happens that the $63$ members of $\B_3$ neatly fall into $7$ classes, each with the same (after adjusting by dividing by $2$)
valuations as the members of $\C_3$. So we can `extract' a letter in $\{1,2,3,4,5,6,7\}$
as well as a member of $\C_3$
in such a way that if you replace the first two columns by that new column,
it is still true that all the rows have the same valuation, in other words, you get a member of $\A_3(n-1)$.
Now keep going, recursively, until we have a  member of $\A_3(1)$, and map its three members to $\{1,2,3\}$.

Of course, there are many choices for these mappings, so we showed how to construct many bijective proofs
to this {\it surprisingly simple} result.
For a table of the bijection for $n$ up to $5$, see the output file
{\tt http://www.math.rutgers.edu/\~{}zeilberg/tokhniot/oBijectionStern1.txt} \quad  .

{\bf Reference}

[S] Richard P. Stanley,
{\it Some linear recurrences motivated by Stern's Diatomic Array}, arXiv:1901.04647v1 [math.CO], 15 January 2019. 
{\tt https://arxiv.org/abs/1901.04647} \quad . \hfill\break
Also in: Amer. Math. Monthly {\bf 127} (2020), 99-111.
%\vfill\eject
\bigskip
\hrule
\bigskip
Shalosh B. Ekhad and Doron Zeilberger, Department of Mathematics, Rutgers University (New Brunswick), Hill Center-Busch Campus, 110 Frelinghuysen
Rd., Piscataway, NJ 08854-8019, USA. \hfill\break
Email: {\tt [ShaloshBEkhad, DoronZeil] at gmail dot com}   \quad .

{\bf Exclusively published in the Personal Journal of Shalosh B. Ekhad and Doron Zeilberger and arxiv.org}

Written: {\bf March 3, 2021}.

\end